\documentclass[12pt]{amsart}
\usepackage{amsmath,amsthm,amsfonts,amssymb,mathrsfs}
\date{\today}
\theoremstyle{definition}
\newtheorem{theorem}{Theorem}
\newtheorem{lemma}{Lemma}
\newtheorem{corollary}{Corollary}
\theoremstyle{definition}
\begin{document}
\title[On countably compact $0$-simple topological inverse semigroups]
{On countably compact $0$-simple topological inverse semigroups}
\footnotetext{This research was supported by the Slovenian
Research Agency grants P1-0292-0101-04 and\break BI-UA/04-06-007. We thank the referee and the editor for comments.}
\author{Oleg~Gutik}
\address{Department of Mathematics, Ivan Franko Lviv National University,
Universytetska 1, Lviv, 79000, Ukraine}
\email{o\_\,gutik@franko.lviv.ua, ovgutik@yahoo.com}
\author{Du\v{s}an~Repov\v{s}}
\address{Institute of Mathematics, Physics and Mechanics, 
and Faculty of Education, University of Ljubljana,
P.O.Box 2964, Ljubljana, 1001, Slovenia}
\email{dusan.repovs@guest.arnes.si}
\begin{abstract}
We describe the structure of $0$-simple countably compact topological inverse
semigroups and the structure of congruence-free countably compact
topological inverse semigroups.
\end{abstract}

\maketitle

We follow the terminology of~\cite{CHK, CP, En}. In this paper all
topological spaces are Hausdorff. If $S$ is a semigroup then we denote
the subset of idempotents of $S$
by
$E(S)$. 
A topological
space $S$ that is algebraically a
semigroup with a continuous
semigroup operation is called a {\em topological semigroup}. A
{\em topological inverse semigroup} is a topological semigroup $S$
that is algebraically an inverse semigroup with continuous
inversion. If $Y$ is a subspace of a topological space $X$ and
$A{\subseteq}Y$, then  we denote by $\operatorname{cl}_Y(A)$
the
topological closure of $A$ in $Y$.

The bicyclic semigroup ${\mathscr{C}}(p,q)$ is the semigroup with
the identity $1$ generated by two elements $p$ and $q$, subject
only to the condition $pq{=}1$. The bicyclic semigroup plays an
important role in the algebraic theory of semigroups and in the
theory of topological semigroups. For example, the well-known
Andersen's result~\cite{Andersen} states that a ($0$--) simple
semigroup is completely ($0$--) simple if and only if it does not
contain the bicyclic semigroup. The bicyclic semigroup admits only
the discrete topology and a topological semigroup $S$ can contain
${\mathscr C}(p,q)$ only as an open subset~\cite{ES}. Neither
stable nor $\Gamma$-compact topological semigroups can contain a
copy of the bicyclic semigroup~\cite{AHK, HildebrantKoch1988}.

Let $S$ be a semigroup and $I_\lambda$  a non-empty set of
cardinality $\lambda$. We define the semigroup operation $' \cdot\, '$ 
on the set $B_{\lambda}(S){=}$
$I_{\lambda}{\times}S^1{\times}I_{\lambda}{\cup}\{ 0\}$ 
as follows
$$
(\alpha, a, \beta){\cdot}(\gamma, b, \delta)=\begin{cases}
(\alpha, ab, \delta), & \text{ if } \beta=\gamma,\\ 0, & \text{ if
} \beta\ne\gamma,
\end{cases}
$$ and $(\alpha, a, \beta){\cdot} 0{=}0{\cdot}(\alpha, a,
\beta){=}0{\cdot} 0{=}0$,
for $\alpha,\beta,\gamma,\delta{\in}
I_\lambda$,
and
$a, b{\in} S^1$. The semigroup $B_{\lambda}(S)$ is
called a~{\it Brandt $\lambda$-extension} of the semigroup
$S$~\cite{GutP23}. Furthermore, if $A{\subseteq}S$ then we shall
denote $A_{\alpha\beta}{=}\{(\alpha, s, \beta)\mid s{\in} A \}$
for $\alpha, \beta{\in} I_{\lambda}$. If a semigroup $S$ is
trivial (i.e. if $S$ contains on\-ly one element), then
$B_\lambda(S)$ is {\it the semigroup of $I_\lambda{\times}
I_\lambda$-matrix units}~\cite{CP}, which we shall denote by
$B_\lambda$. By Theorem~3.9 of \cite{CP}, an inverse semigroup $T$ is
completely $0$-simple if and only if $T$ is isomorphic to a Brandt
$\lambda$-extension $B_{\lambda}(G)$ of some group $G$ and
$\lambda{\geqslant}1$. We also note
that if $\lambda{=}1$, then
the semigroup $B_{\lambda}(S)$ is isomorphic to the semigroup $S$
with adjoint zero. Gutik and Pavlyk \cite{GutP68}
proved that
any continuous homomorphism from the infinite topological
semigroup of matrix units into a compact topological semigroup is
annihilating, and hence the infinite topological semigroup of
matrix units does not embed into a compact topological semigroup.
They also
showed that if a topological inverse semigroup $S$
contains a semigroup of matrix units $B_{\lambda}$, then
$B_{\lambda}$ is a closed subsemigroup of~$S$.

Suschkewitsch \cite{Suschkewitsch1928}  proved that any finite
semigroup $S$ contains a minimal ideal $K$. He also showed that
$K$ is a completely simple semigroup and described the structure
of finite simple semigroups. Rees \cite{Rees1940}  generalized
the Suschkewitsch Theorem and showed that if a semigroup $S$
contains a minimal ideal $K$ then $K$ is isomorphic to a Rees
matrix semigroup $M[G;I,\Lambda,P]$ over a group $G$ with a
regular sandwich matrix $P$. He also proved that any completely
$0$-simple semigroup is isomorphic to a Rees matrix semigroup
$M[G;I,\Lambda,P]$ over a $0$-group $G^0$ with a regular sandwich
matrix $P$. Wallace~\cite{Wallace1956} proved the topological
analogue of
the Suschkewitsch-Rees Theorem for compact topological
semigroups: \emph{every compact topological semigroup contains a
minimal ideal, which 
is topologically isomorphic to a
topological paragroup}. Paalman-de-Miranda
\cite{Paalman-de-Miranda1964} proved that any $0$-simple compact
to\-po\-logical semigroup $S$ is completely $0$-simple, the zero
of $S$ is an isolated point in $S$ and $S{\setminus}\{ 0\}$ is
homeomorphic to the topological product $X{\times}G{\times}Y$,
where $X$ and $Y$ are compact topological spaces and $G$ is
homeomorphic to the underlying space of a maximal subgroup of $S$,
contained in $S{\setminus}\{ 0\}$. Owen \cite{Owen} showed that
if $S$ a locally compact completely simple topological semigroup,
then $S$ has a structure similar
to a compact simple topological
semigroup. Owen also gave an example which shows that a similar
statement does not hold for a locally compact completely
$0$-simple topological semigroup.  Gutik and
Pav\-lyk \cite{GutP68}
proved that the subsemigroup of idempotents of a compact
$0$-simple topological inverse semigroup is finite, and hence the
topological space of a compact $0$-simple topological inver\-se
semigroup is homeomorphic to a finite topological sum of compact
topological group and a single point.

A Hausdorff topological space $X$ is called \emph{countably
compact} if any open countable cover of $X$ contains a finite
subcover~\cite{En}. In this paper we shall prove that the bicyclic
semigroup cannot be embedded into any countably compact
topological inverse semigroup. We shall also describe the
structure of $0$-simple countably compact topological inverse
semigroups and the structure of congruence-free countably compact
topological inverse semigroups.

\begin{theorem}\label{th1.3}
{\it A countably compact topological inverse semigroup cannot contain
the bicyclic semigroup. Therefore every (0-)simple countably
compact topological inverse semigroup is (0-)completely simple.}
\end{theorem}

\begin{proof}
Let $T$ be a countably compact topological inverse semigroup and
suppose that $T$ contains ${\mathscr{C}}(p,q)$ as a subsemigroup.
Let $S=\operatorname{cl}_T({\mathscr{C}}(p,q))$. Then by
Theorem~3.10.4 of 
\cite{En}, $S$ is a countably compact space and by
Proposition~II.2 of \cite{ES}, $S$ is a topological inverse semigroup.
Thus by Corollary~I.2 of \cite{ES}, the semigroup ${\mathscr{C}}(p,q)$
is a discrete subspace of $S$ and by Theorem~I.3 of \cite{ES},
${\mathscr{C}}(p,q)$ is an open subspace of $S$ and
$S{\setminus}{\mathscr{C}}(p,q)$ is an ideal in $S$. Therefore any
element of ${\mathscr{C}}(p,q)$ is an isolated point in the
topological space $S$. We define the maps $\varphi{\colon} S{\to}
E(S)$ and $\psi{\colon} S{\to} E(S)$ by the formulae
$\varphi(x){=}xx^{-1}$ and $\psi(x){=}x^{-1}x$. Since
$S{\setminus}{\mathscr{C}}(p,q)$ is an ideal of $S$,
$A{=}\varphi^{-1}(\{ 1\})\cup\psi^{-1}(\{
1\}){\subseteq}{\mathscr{C}}(p,q)$, and since the maps $\varphi$
and $\psi$ are continuous $A$ is a clopen and hence countably
compact infinite subset of $S$. But $A$ is an open subspace of $S$
whose elements are isolated points in $S$. A contradiction.

The second part of the theorem follows from
Theorem~2.54 of \cite{CP}.
\end{proof}

Let $\mathscr{S}$ be a class of topological semigroups. Let
$\lambda$ be a cardinal ${\geqslant}1$, and
$(S,\tau){\in}\mathscr{S}$. Let $\tau_{B}$ be a topology on
$B_{\lambda}(S)$ such that $\left(B_{\lambda}(S),
\tau_{B}\right){\in}\mathscr{S}$ and $\tau_{B}|_{(\alpha, S,
\alpha)}{=}\tau$ for some $\alpha\in I_{\lambda}$. Then
$\left(B_{\lambda}(S), \tau_{B}\right)$ is called a {\it
topological Brandt $\lambda$-extension of $(S, \tau)$} in
$\mathscr{S}$ \cite{GutP23}.

Let $\alpha,\beta,\gamma,\delta{\in} I_{\lambda}$ and $A$ be a
subspace of $S$. Since the restriction
$\varphi_{\alpha\beta}^{\gamma\delta}\big|_{A_{\alpha\beta}}\colon
A_{\alpha\beta}\to A_{\gamma\delta}$ of the map
$\varphi_{\alpha\beta}^{\gamma\delta}\colon B_{\lambda}(S){\to}
B_{\lambda}(S)$ defined by the formula
$\varphi_{\alpha\beta}^{\gamma\delta}(s){=}(\gamma,1,\alpha)\cdot
s\cdot(\beta,1,\delta)$ is a homeomorphism, we get the following:

\begin{lemma}\label{lem} 
{\it Let $\lambda{\geqslant} 1$ and $B_{\lambda}(S)$ be a
topological Brandt $\lambda$-extension of a topological semigroup
$S$ and $A$  a subspace of $S$. Then the subspaces
$A_{\alpha\beta}$ and $A_{\gamma\delta}$ in $B_{\lambda}(S)$ are
homeomorphic for all $\alpha,\beta,\gamma,\delta{\in}
I_{\lambda}$.}
\end{lemma}

\begin{theorem}\label{th1.5}
{\it Let $S$ be a $0$-simple countably compact topological inverse
semigroup. Then there exist a nonempty finite set $I_{\lambda}$ of
cardinality $\lambda$ and a countably compact topological group
$H$ such that $S$ is topologically isomorphic to a topological
Brandt $\lambda$-extension $B_{\lambda}(H)$ of $H$ in the class of
topological inverse semigroups. Moreover, $S$ is homeomorphic to a
finite topological sum of countable compact topological groups and
a single point.}
\end{theorem}

\begin{proof}
By Theorem~\ref{th1.3}, the semigroup $S$ is completely $0$-simple.
Now Theorem~3.9 of \cite{CP} implies that there exist a nonempty set
$I_{\lambda}$ of cardinality $\lambda$ and a group $G$ such that
$S$ is algebraically isomorphic to $B_{\lambda}(G)$. Therefore for
any $\alpha{\in} I_{\lambda}$ the subset $G_{\alpha\alpha}$ is a
subgroup of $B_{\lambda}(G)$ and since $B_{\lambda}(G)$ is a
topological inverse semigroup, a topological subspace
$G_{\alpha\alpha}$ of $B_{\lambda}(G)$ with the induced
multiplication is a topological group. We fix $\alpha{\in}
I_{\lambda}$ an put $H{=}G_{\alpha\alpha}$. Then the
topo\-lo\-gical semigroup $S$ is topologically isomorphic to a
topological Brandt $\lambda$-extension $B_{\lambda}(H)$ of the
topological group $H$.

Let $e_H$ be the identity of $H$. Then the subsemigroup
 $
B_{\lambda}(e_{H}){=}\{
0\}\cup\{(\alpha,e_{H},\beta)\mid\alpha,\beta{\in} I_{\lambda}\}
 $
of $B_{\lambda}(H)$ is algebraically isomorphic to the semigroup
of matrix units $B_{\lambda}$. By Theorem~14 \cite{GutP68},
$B_{\lambda}(e_{H})$ is a closed subsemigroup of $B_{\lambda}(H)$
and hence by Theorem~3.10.4 of \cite{En}, $B_{\lambda}(e_{H})$ is a
countably compact topological space. Therefore
Theorem~6 of \cite{GutP68} implies that $B_{\lambda}(e_{H})$ is a
finite discrete subsemigroup of $B_{\lambda}(H)$ and hence the set
$I_{\lambda}$ is finite.

We define the maps $\varphi\colon B_{\lambda}(H){\to}
B_{\lambda}(e_{H})$ and $\psi\colon B_{\lambda}(H){\to}
B_{\lambda}(e_{H})$ by the formulae $\varphi(x)=xx^{-1}$ and
$\psi(x){=}x^{-1}x$. Since $B_{\lambda}(H)$ is a topological
inverse semigroup the maps $\varphi$ and $\psi$ continuous and
hence by Lemma~4 of \cite{GutP68}, the set $H_{\alpha\beta}{=}
\varphi^{-1}((\alpha,e_{H},\beta)){\cap}\varphi^{-1}((\alpha,e_{H},\beta))$
is clo\-pen in $B_{\lambda}(H)$. By Lemma~\ref{lem}, the subspaces
$H_{\alpha\beta}$ and $H_{\gamma\delta}$ are homeomorphic for any
$\alpha,\beta,\gamma,\delta{\in} I_{\lambda}$, and hence all of
them are homeomorphic to the topological group $H$.
\end{proof}

A Tychonoff topological space $X$ is called \emph{pseudocompact}
if every continuous real-valued function on $X$ is bounded. Since
the topological space of $T_0$-topological group is Tychonoff and
any topological sum of Tychonoff spaces is a Tychonoff 
space, Theorem~3.10.20 of \cite{En} implies:

\begin{corollary}\label{cor1.8}
{\it The topological space of a $0$-simple countably compact
topological inverse semigroup is Tychonoff and hence
pseudocompact.}
\end{corollary}

Let $X$ be a topological space. The pair $(Y, c)$, where $Y$ is a
compactum and $c{\colon} X{\to} X$ is a homeomorphic embedding of
$X$ into $Y$, such that $\operatorname{cl}_{Y}c(X){=}Y$, is called
a \emph{compactification} of the space $X$. Define the ordering
$\preccurlyeq$ on the family ${\mathcal{C}}(X)$ of all
compactifications of a topological space $X$ as
follows: $
c_2(X){\preccurlyeq} c_1(X)$ if and only if there exists a
continuous map $f{\colon} c_1(X){\to} c_2(X)$ such that
$fc_1{=}c_2$. The greatest element of the family
${\mathcal{C}}(X)$ with respect to the ordering $\preccurlyeq$ is
called the \emph{Stone-\v{C}ech compactification} of the space $X$
and it is denoted by $\beta{X}$.  Comfort and
Ross \cite{ComfortRoss}
proved that the Stone-\v{C}ech compactification of a
pseudocompact topological group is a topological group. The next
theorem is an analogue of the Comfort--Ross Theorem:

\begin{theorem}\label{th1.9}
{\it Let $S$ be a $0$-simple countable compact topological inverse
semigroup. Then the Stone-\v{C}ech compactification of $S$ admits
a structure of $0$-simple topological inverse semigroup with respect
to which the inclusion mapping of $S$ into $\beta{S}$ is a
topological isomorphism.}
\end{theorem}

\begin{proof}
By Theorem~\ref{th1.5}, $S$ is topologically isomorphic to a Brandt
$\lambda$-extension of some topological group $H$ in the class of
topological inverse semigroups and $\lambda{<}\omega$. Now by
Lemma~\ref{lem}, the subspaces $H_{\alpha\beta}$ and
$H_{\gamma\delta}$ are homeomorphic in $B_{\lambda}(H)$, for any
$\alpha,\beta,\gamma,\delta{\in} I_{\lambda}$. Since a maximal
subgroup in $S$ is closed we have that $H_{\alpha\beta}$ is a
clopen subset of $B_{\lambda}(H)$, for every $\alpha,\beta{\in}
I_{\lambda}$. By Corollary~\ref{cor1.8}, the topological space
$B_{\lambda}(H)$ is pseudocompact. Since any clopen subspace of a
pseudocompact topological space is pseudocompact
(see ~\cite{Colmex}) the subspace $H_{\alpha\beta}$ is
pseudocompact, for every $\alpha,\beta{\in} I_{\lambda}$.
Obviously, the topological space $B_{\lambda}(H){\setminus}\{ 0\}$
is homeomorphic to $H{\times} I_{\lambda}{\times} I_{\lambda}$.
Since the topological space $I_{\lambda}{\times} I_{\lambda}$ is
finite and hence compact, by Corollary~3.10.27 of \cite{En}, the space
$B_{\lambda}(H){\setminus}\{ 0\}$ is pseudocompact. Now by
Theorem~1 of \cite{Glicksberg}, we have
 $
\beta(H{\times} I_{\lambda}{\times} I_{\lambda}){=}\beta
H{\times}\beta I_{\lambda}{\times}\beta I_{\lambda}{=}\beta
H{\times} I_{\lambda}{\times} I_{\lambda}
 $
and therefore $\beta(B_{\lambda}(H)){=}B_{\lambda}(\beta H)$.
\end{proof}

\begin{corollary}\label{cor1.10}
{\it Every $0$-simple countable compact topological inverse semigroup
is a dense subsemigroup of a $0$-simple compact topological
inverse semigroup.}
\end{corollary}

If $S$ is completely simple inverse semigroup then the semigroup
$S$ with joined zero $S^0$ is completely $0$-simple and hence by
Theorem~3.9 of \cite{CP}, the semigroup $S^0$ is isomorphic to a
Brandt $\lambda$-extension $B_{\lambda}(G)$ of some group $G$.
Therefore any nonzero idempotent of $S^0$ is primitive. Let $e$
and $f$ are nonzero idempotents of $S^0$. Since $S$ is an inverse
subsemigroup of $S^0$ we have $ef{=}fe{\leqslant} e$ and
$ef{=}fe{\leqslant} f$, and hence $e{=}ef{=}f$. Thus, the inverse
semigroup $S$ contains the unique idempotent and hence it is a
group. Therefore a completely simple inverse semigroup is a group
and Theorem~\ref{th1.3} implies that \emph{every simple countable
compact topological inverse semigroup is a topological group}.

A~semigroup $S$ is called {\it congruence-free} if it has only two
congruences: the identity relation and the universal
relation~\cite{Sch}.

\begin{theorem}\label{th1.12}
{\it Let $S$ be a congruence-free countably compact topological inverse
semigroup with zero. Then $S$ is isomorphic to a finite semigroup
of matrix units.}
\end{theorem}

\begin{proof} Suppose not. Since the semigroup $S$ contains a
zero by Theorem~\ref{th1.5}, $S$ is topologically isomorphic to a
topological Brandt $\lambda$-extension $B_{\lambda}(H)$ of a
pseudocompact topological group $H$ in the class of topological
inverse semigroups and $\lambda{<}\omega$. Suppose that the
group $H$ is not  trivial. Then we define a map $h\colon
B_{\lambda}(H){\rightarrow} B_{\lambda}$ by the formulae
$h((\alpha,g,\beta)){=}(\alpha,\beta)$ and $h(0){=}0$. Since
 $
h((\alpha,g,\beta)(\gamma,s,\delta)){=h}((\alpha,gs,\delta)){=}
(\alpha,\delta)=(\alpha,\beta)(\gamma,\delta){=}$ $
h((\alpha,g,\beta))h((\gamma,s,\delta))
 $
for $\beta{=}\gamma$ and
 $
h((\alpha,g,\beta)(\gamma,s,\delta))=h(0){=} 0{=}
(\alpha,\beta)(\gamma,\delta)=
h((\alpha,g,\beta))h((\gamma,s,\delta))
 $
for $\beta\neq\gamma$, the map $h$ is a homomorphism. This
contradicts the assumption that $S$ is a congruence-free
semigroup.
\end{proof}

\end{document}